# Challenges and perspectives in LMD: comparative study and intelligent mobility proposal

David García-Retuerta

**Abstract**

The "Last Mile Delivery" (LMD) refers to the last and most inefficient part of the supply chain. This is caused by the spatial distribution of disperse small receiving points, the ever-growing demand for faster shipment and the new time constraints of deliveries. Moreover, the small urban vehicles used for package distribution make LMD the most polluting part of the supply chain. This study describes the existing methods to improve the efficiency of LMD, its challenges and makes a conceptual proposal. The state-of-the-art techniques are categorised according to their data sources, the theoretical framework of the proposed models and the dynamic nature of the solutions.

## 1. Introduction

Home delivery of goods in urban areas is becoming a growing problem, which is considered of vital importance for the user experience in online shipping (Xiao, Wang, Lenzer, & Sun, 2017). The movement of goods using trains or ships is agreed to be cost-efficient currently (Wang, Zhang, Liu, Shen, & Lee, 2016). However, when the goods arrive at the destination port or warehouses, they must be distributed to their recipients. This last step in the freight supply chain is called Last Mile Delivery (LMD) and is considered to be the most inefficient, polluting and cost-inefficient part (Moroza & Polkowski, 2016) of the supply chain. An integrated overview of last mile transport can be found in (Juhász & Bányai, 2018). The cost of last-mile shipping is estimated to range from 13% to 75% of total shipping costs (Gevaers, Voorde, & Vanelslander, 2009), and can account for up to 40% of the price paid by the buyer (Kong, Li, Yu, Zhao, & Huang, 2017), leading to interest in increasing its efficiency. The latter is perhaps especially important in urban areas, where traffic conditions are unstable, and the unpredictable behaviour of stakeholders (drivers, cyclists and pedestrians) can make it difficult to meet delivery schedules.

The problems caused by sending small packages -the majority of packages sold by Amazon weigh less than 2.2 kg. (Chen & Chankov, 2017)- is likely to increase in the next few years due to the increasing relevance of e-commerce. Reviews of the impact of home delivery on urban freight transport can be found in (Savelsbergh & Woensel, 2016) and also in (Visser, Nemoto, & Browne, 2014). Moreover, sellers have to cope with buyers' increased demands on the level of service: shorter delivery times, same-day shipments and shipments on demand, which can reduce the delivery time to one or two hours (Arslan, Agatz, Kroon, & Zuidwijk, 2016).

This level of demand is a problem for both delivery companies and city residents who suffer its negative effects, such as: a higher price of customer service, traffic congestion, safety issues and increased environmental pollution (Moroza & Polkowski, 2016) (Laranjeiro, et al., 2019) (Perboli & Rosano, 2019) (Ranieri, Digiesi, Silvestri, & Roccotelli, 2018).

From an economic point of view, the need to increase efficiency due to decreasing marginal revenues has led several stakeholders to consider alternatives to increase efficiency in last mile delivery. Among the most classical alternatives is optimisation in different forms: optimisation of vehicle routes and combining the transport of goods with the transport of people, for example by bus (Pimentel & Alvelos, 2018) or taxi (Chen & Pan, 2016) (Li B. , Krushinsky, Woensel, & Reijers, 2016). A crowdsourcing approach can also be attempted, through occasional drivers (via crowdsourcing or in-store customers) (Gdowska, Viana, & Pedroso, 2018), (Arslan A. , Agatz, Kroon, & Zuidwijk, 2016), (Devari, Nikolaev, & He, 2017) (Chen, Cheng, Lau, & Misra, 2015). In (Behrend & Meisel, 2018) the increase in profit and service level obtained by using a crowdsourced delivery platform is calculated. This follows the current trend of developing the collaborative economy in big cities (Cohen, B., Almirall, E., & Chesbrough, H., 2016).

Other alternatives proposed include, for example, (Perboli & Rosano, 2019) (Janjevic, Winkenbach, & Merchán, 2019) the use of parcel lockers – secure and unattended, with independent collection - or attended collection points - with own or third-party premises - as well as the use of alternative transport methods such as bicycles, electric vehicles, autonomous vehicles or drones.

Another recent consideration is related to environmental impact, as the delivery of goods in urban areas is considered to be responsible for 25% of $CO_2$ emissions and 30-50% of other pollutants (ALICE / ERTRAC Urban mobility WG, 2015). Similar values are mentioned in (Muñoz-Villamizar, Montoya-Torres, & Vega-Mejía, 2015).

The LMD challenge is intrinsically equivalent to the well-known *Travelling salesman problem*, which fundamentally tried to answer the question: "Given a list of cities and the distances between each pair of cities, what is the shortest possible route that visits each city exactly once and returns to the origin city?" (G. Reinelt, 2003) (M. Yu, 2019). Hence, all the knowledge already gathered in this classical research topic can be used in LMD, benefiting from classical solutions.

Our work diverges from the current literature as it proposes a categorization based on new classification actions, which are then used to put forward a theorical LMD proposal. The studied works are organised according to the following categories: (1) data sources, either real or synthetical   (2) algorithmic foundation, either heuristic, mathematical, or a mixture of both   (3) static or dynamic nature of the algorithm. Moreover, an overall categorization is performed, and urban road freight impact is analysed in a dedicated section due to its high importance.

The article is organised as follows: Section 2 describes the most common challenges found by researchers when solving the last mile delivery problem. Section 3 carries out a comparative study of the most relevant literature, analysing the studies based on their data source and on the used technical solution. Section 4 generalizes the previous classification and summarizes the most important characteristics of the analysed research works. Finally, a LMD solution is proposed in Section 5, and in Section 7 we discuss the findings and thoughts this article provides.

## 2. Challenges

A good description of the last mile delivery problem can be found in (Souza, Goh, Lau, Ng, & Tan, 2014). In their study, they describe the last mile delivery challenges in the city-state of Singapore and also outline a possible solution. According to the authors, the set of issues listed for the LMD problem may include:

- Pick-up/delivery points are often located in restricted access areas away from major distribution centres, dispersed in clusters of disparate demand or supply and subject to regulated access (e.g., slots).
- Congestion.
- Fleet/cargo usage limits.
- Dynamic interaction between many competing interests and services, policies and interventions.

The problem described for Singapore, which is shared with other cities, includes the existence of several competing delivery companies (logistics service providers) in the city. Thus, the resulting last mile logistics model is fragmented and uncoordinated, which leads to:

(i) low truck capacity utilisation.

(ii) excessive truck movements.

(iii) higher cost for the whole system.

(iv) negative environmental impacts.

The authors indicate that there are some benefits to be gained from having Urban Consolidation Centres (UCCs), such as:

(a) increased load factor.

(b) decrease in commercial vehicle traffic.

(c) reduced emissions.

(d) increased levels of service.

However, according to the authors, the future of UCCs is not guaranteed without the support of a governmental authority due to their voluntary nature. Based on the idea of a UCC, the authors propose a solution focused on cooperation, while optimising the resources of the respective companies and also strengthening their own market position. They propose to establish a UCC as a (virtual) marketplace, pooling requests and generating viable vehicle routing options (both efficient and sustainable) within existing lead times and constraints. In the proposed scenario, relevant data should be harmonised and exchanged to create a common plan to increase truck capacity utilisation and asset substitution capacity. This electronic market is proposed to function as a combined iterative auction mechanism in which bid generation takes place within each agent, and a centralised auctioneer carries out a price adjustment procedure.

The authors plan four areas of interest for achieving collaboration in urban logistics (Souza, Goh, Lau, Ng, & Tan, 2014):

- green collaborative delivery.
- multi-objective synchronisation and planning.
- multi-stakeholder coordination.
- data harmonisation and analysis.

The authors mention some variants of the Vehicle Routing Problem (VRP) to obtain a solution based on mathematical models, such as the VRP with time windows, the Capacitated VRP, the multi-depot VRP, the dial-a-ride VRP, the heterogeneous fleet VRP, the VRP with pickup and delivery, among other. However, as it falls within the combinatorics field, with a vast number of possible solutions to be considered, they are all discarded.

According to the authors, the VRP problem should be extended to address environmental impacts and an emission estimation model should be established. To obtain such a model, the possibility of dynamic routing and scheduling for a vehicle based on real-time traffic conditions for a route with multiple pick-up and drop-off points is suggested.

## 3. Comparative Study of the literature

The carried-out review shows that the literature of the last mile delivery is diverse and rapidly emerging, with new authors proposing innovative approaches and incremental improvements of existing solutions.

In addition, the subject is also rather fragmented, so our aim in the current study is to create some order to it, grouping the works by their theoretical foundation into: heuristic-based, mathematically-based and a mixture of the previous. This classification was chosen to provide researchers a clear image of the advantages and disadvantages of making use of approximation functions, either fully or partially. As a result, researchers can choose the method which better suits their considered situation and use it in order to solve their problem efficiently.

Moreover, we make a distinction between static and dynamic solutions, and propose an overall in-detail classification. This separation focuses on the random nature of the situation: if the considered problem is unstable, a dynamic solution will better adapt to the new states; but if the problem is stable, it would only increase its complexity.

Data sources are also considered, as either use real-world data or synthetical data, depending on availability and the development timeline. The current section explores all the former points.

### 3.1. Common Data Sources

An interesting reflection on the data used to test the different proposals can be found in (Perboli, Rosano, Saint-Guillain, & Rizzo, 2018), where they lack a factual dataset to make realistic comparisons for the VRP problem. They find that in realistic (urban) scenarios, the quality of the solution depends on several parameters involved:

1. the geographical distribution of the target groups.

2. the characteristics of the different types of vehicles.

3. the use of parcel lockers to deliver part of the demand.

Overall, a considerable number of studies use synthetically generated data. They can be divided into two main groups:

- Studies using common or standard data, which can help establish a basis for a comparison of different proposals, which eliminates variability in the data.

- Studies that generate their own data.

| | Sao Paulo (BR) | 15500 individual shipments daily | Hypothetical topology | Multiple shipping companies | UCC / multilevel networks | Simulation | Environmental considerations | Yokohama (JP) | Berlin (DE) | Food delivery | Parcel delivery | Turin (IT) | 100, 250 and 500 destinations | 15 mails + 5 greens | Monte Carlo | San Francisco (USA) | Washington DC (USA), UPS locations | Bremen (DE) map only | Bogotá (CO) | 61 stores, 3 different companies | Vienna (AT) | e-groceries, 255 stores, 24 vehicles | Rotterdam (NL) | Synthetic data |
|---|---|---|---|---|---|---|---|---|---|---|---|---|---|---|---|---|---|---|---|---|---|---|---|---|
| (Janjevic, Winkenbach, & Merchán, 2019) | × | × | × | | | | | | | | | | | | | | | | | | | | | |
| (Firdausiyah, Taniguchi, & Qureshi, 2019) | | | | × | × | × | × | × | | | | | | | | | | | | | | | | × |
| (Chatterjee, Greulich, & Edelkamp, 2016) | | | | | | × | × | | | | | | | | × | | | × | | | | | | × |
| (Fikar, 2018) | | | | | | | | | | | | | | | | | | | | | × | × | | × |
| (Martins-Turner & Nagel, 2019) | | | | × | | | × | | | × | × | | | | | | | | | | | | | |
| (Rosano, Demartini, Lamberti, & Perboli, 2018) | | | | | | × | × | | | | | × | × | × | × | | | | | | | | | × |
| (Li B. , Krushinsky, Woensel, & Reijers, 2016) | | | | | | | × | | | | | | | | | × | | | | | | | | |
| (Chen & Chankov, 2017) | | | | | | | | | | | | | | | | | | × | | | | | | × |
| (Muñoz-Villamizar, Montoya-Torres, & Vega-Mejía, 2015) | | | | | | | × | | | | | | | | | | | | × | × | | | | × |

includes information on validations using real data or synthetic data generated from real data or configurations.

| | Sao Paulo (BR) | 15500 individual shipments daily | Hypothetical topology | Multiple shipping companies | UCC / multilevel networks | Simulation | Environmental considerations | Yokohama (JP) | Berlin (DE) | Food delivery | Parcel delivery | Turin (IT) | 100, 250 and 500 destinations | 15 mails + 5 greens | Monte Carlo | San Francisco (USA) | Washington DC (USA), UPS locations | Bremen (DE) map only | Bogotá (CO) | 61 stores, 3 different companies | Vienna (AT) | e-groceries, 255 stores, 24 vehicles | Rotterdam (NL) | Synthetic data |
|---|---|---|---|---|---|---|---|---|---|---|---|---|---|---|---|---|---|---|---|---|---|---|---|---|
| (Janjevic, Winkenbach, & Merchán, 2019) | ✕ | ✕ | ✕ | | | | | | | | | | | | | | | | | | | | | |
| (Firdausiyah, Taniguchi, & Qureshi, 2019) | | | ✕ | ✕ | ✕ | ✕ | ✕ | ✕ | | | | | | | | | | | | | | | | ✕ |
| (Chatterjee, Greulich, & Edelkamp, 2016) | | | | | ✕ | ✕ | | | | | | | | | ✕ | | | ✕ | | | | | | ✕ |
| (Fikar, 2018) | | | | | | | | | | | | | | | | | | | | | ✕ | ✕ | | ✕ |
| (Martins-Turner & Nagel, 2019) | | | | ✕ | | ✕ | | | | ✕ | ✕ | | | | | | | | | | | | | |
| (Rosano, Demartini, Lamberti, & Perboli, 2018) | | | | | ✕ | ✕ | | | | | | ✕ | ✕ | ✕ | ✕ | | | | | | | | | ✕ |
| (Li B. , Krushinsky, Woensel, & Reijers, 2016) | | | | | | ✕ | | | | | | | | | | ✕ | | | | | | | | |
| (Chen & Chankov, 2017) | | | | | | | | | | | | | | | | | ✕ | | | | | | | ✕ |
| (Muñoz-Villamizar, Montoya-Torres, & Vega-Mejía, 2015) | | | | | | ✕ | | | | | | | | | | | | | ✕ | ✕ | | | | ✕ |
| (Arslan A. , Agatz, Kroon, & Zuidwijk, 2016) | | | | | | | | | | | | | | | | | | | | | | | ✕ | ✕ |
| (Jiang, L., Zang, X, *et al.*, 2022) | | | ✕ | | | ✕ | | | | | ✕ | | | | | | | | | | | | | ✕ |

**Table 1. Validations performed with synthetic and real-world data. Proposals in the vertical axis and carried out validations in the horizontal axis**

### 3.1.1. Studies Based on Real-World Data

The literature presents many studies related to the LMD problem, where important players of the logistics business share their experiences in several parts of the world. Some of the most prominent examples, which provide the most useful insights are listed below.

(Janjevic, Winkenbach, & Merchán, 2019) presents a study based on data collected from B2W, the largest last-mile distributor in São Paulo, Brazil. Their shipment vehicles are vans with a capacity of 2.4 m$^3$. The authors identify how many satellites and when Collection & Delivery Points (CDPs) are needed for their network. This calculation is based on data from other studies, which indicate that 10% of Brazilians prefer to go to pick up deliveries rather than wait for them at home (5% prefer to use lockers and another 5% prefer CDPs). To assess the position that CDPs should have, 99 configurations are tested. The authors estimate the savings introduced by CDPs (compared to not having them) to be between 4.6 and 5.0% of the cost. In addition to this analysis, a sensitivity analysis is performed to study the feasibility of integrating CDPs into the network: varying parameters such as the percentage of users who prefer to go to the CDP for delivery, or the number of active CDPs; the authors found that the overall cost (and therefore the profitability of the model) is very sensitive to these variations.

In (Zhou, Baldacci, Vigo, & Wang, 2018) they also conduct a performance analysis of their proposal, based on a genetic algorithm, using both real and randomly generated data. In the city of Chongqing (China), two delivery companies form an alliance. A network of satellites and Customer Pickup (CP) points - a multilevel model - is designed with: two warehouses, 12 satellites, 40 CPs and a total of 164 recipients. The authors calculate the results in 4 different scenarios:

1. Joint distribution with two delivery options (home or CP).
2. Joint distribution with no delivery option (always home delivery).
3. Independent distribution with delivery options,
4. Independent distribution without delivery options.

The results of the experiment have been summarised in Table 2, where it can be seen that the joint distribution allows a reduction of 7.1 and 7.6% respectively. In addition, allowing the option of collection at the CP results in a cost reduction of about 9%, and the number of vehicles is considerably reduced (from 37 to 20 in the best case). This study is complemented by a sensitivity analysis for the different parameters.

| Scenario | Total cost | Connection cost | Satellites number | CP number | Number of 1 N routes | Number of 2 N routes |
|---|---|---|---|---|---|---|
| 1 | 9395.0 | 2506.6 | 6 | 27 | 3 | 20 |
| 2 | 10317.1 | 0 | 6 | 0 | 3 | 37 |

| 3 | 10113.5 | 1089.7 | 6 | 18 | 3 | 30 |
| 4 | 11166.8 | 0 | 7 | 0 | 3 | 39 |

**Table 2. Results of the experiment with real data. For each of the 4 scenarios defined, the total cost, connection cost, number of satellites, number of collection points, number of first level routes, number of second level routes and number of connections are provided.**

In addition, the authors perform a second study is performed based on a data generator for multi-level networks, where facilities and destinations are created within concentric circles of different radius. The depots are located in area 1, the satellites in area 2 or 3 and the CPs are placed in area 3, and the destinations in area 4 around the CPs. All locations are randomly generated (Figure 1). The synthetic dataset is used to test the impact of the different elements of the proposed solution. This study is also complemented by a sensitivity analysis for the different parameters.

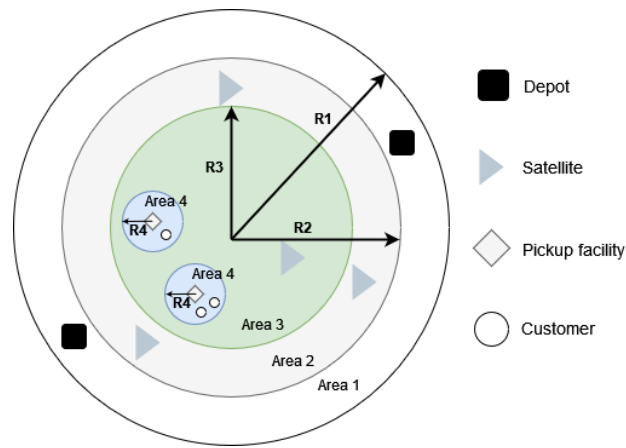

**Figure 1. Schematic representation for a multi-level network to generate data.**

(Martins-Turner & Nagel, 2019) analyses two scenarios to study the performance of their proposal. The <u>first scenario</u> is based on data on food distribution to retailers in the city of Berlin (DE), with 1928 requests, 15 different companies, 867 shops and 17 distribution centres. Food distribution has some special characteristics, e.g., different types of food (fresh, frozen, dry) may need different types of transport, which implies in this case 15 x 3 different delivery companies. Vehicles of one type cannot be reused for shipments of another type.

The effect of being able to make only one or several trips is studied, and time restrictions for delivery are established. Several configurations of the dispatch network are tested, including access tolls for low-emission zones, only electric vehicles or a UCC for the low-emission zone. The results obtained are: being able to make more than one trip reduces the vehicle fleet by 62.5-68% and the cost decreases by 35-44%, while the distance and travel time increases by 6%.

The underlined second scenario is based on parcel delivery in one particular district of Berlin: Wilmersdorf. Two variants are also introduced (Figure 2):

1. Direct distribution from the warehouse using light vehicles.
2. Distribution to intermediate micro warehouses using light vehicles and from these to the destination using cargo bikes (some parcels are still shipped directly).

In this case, the authors conclude that being able to make several trips makes it possible to greatly reduce the number of vehicles required and therefore the cost, due to the low carrying capacity of the bicycles. The use of light trucks is also reduced by supplying micro warehouses and delivering only private parcels.

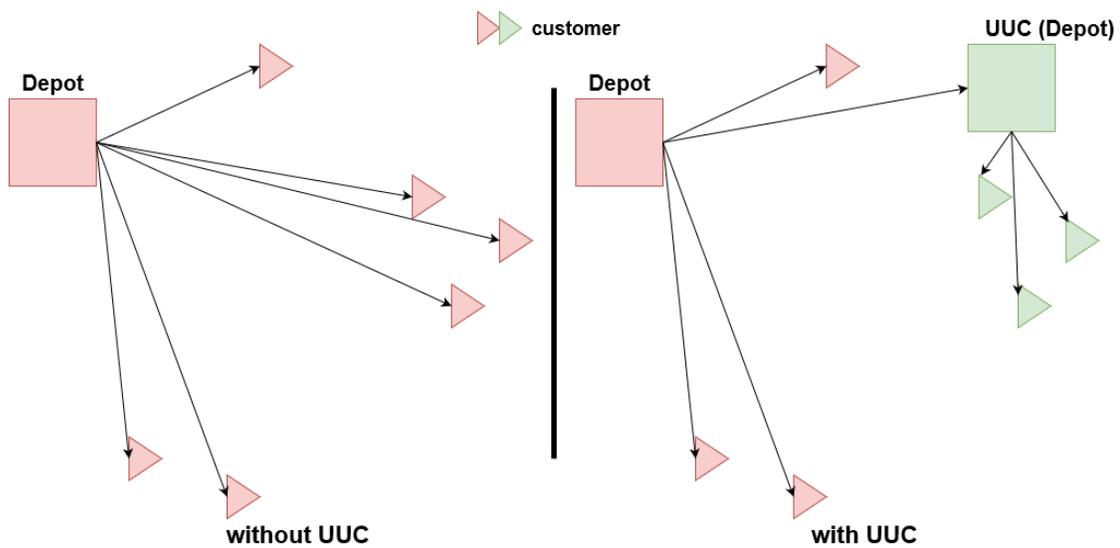

**Figure 2. Transport without and with UCC between the warehouse and the customer.**

(Rosano, Demartini, Lamberti, & Perboli, 2018) presents a proposal for a platform to connect customers with parcel delivery companies (pickup & delivery). The proposal is derived from the real database on the distribution of destinations and volumes, provided by an international delivery company in the city of Turin (Italy). For testing purposes, 360 datasets of three different sizes (500, 250 and 100 destinations) are randomly generated. There are 15 delivery companies and 5 green delivery companies (using cargo bikes) operating in Turin. The platform considers the dynamics of accepting new requests once routes are already started. The route optimisation process involves solving several problems: Vehicle Routing Problem with Time Windows (VRPTW) combined with load balancing, a stochastic Traveling Salesman Problem (TSP) and a dynamic stochastic VRPTW. The result includes the calculation of several KPIs (Key Performance Indicators) to assess the sustainability of the proposal. With a stop-to-stop distance of about 700 metres and a minimum of 80 dispatches per day, the operating cost per stop is calculated to decrease by 25%, the social cost of emissions to decrease by 21% and the number of dispatches per hour to increase by 37%.

(Li B., Krushinsky, Woensel, & Reijers, 2016) propose a parcel delivery solution though the model Share-a-Ride Problem (SARP) using taxis, based on data collected from real-world taxis of the city of San Francisco (USA). They introduce two different variations in their proposed model:

1. SARPs with stochastic travel times, due to traffic conditions.
2. SARP with stochastic destination, which may occur in some areas, where no precise destination address is known but an approximate one (rural areas in some countries, for example).

The stochastic problem is solved in two steps: first, the deterministic part of the problem (the vehicle route) is solved and then, the random part (the travel times or the precise destinations) is incorporated via scenario simulations. Once the random variables are known, it is possible to check whether the constraints for the times (both delivery and travel times for the taxi occupants) are satisfied. The following elements are incorporated in the cost function:

1. Penalties for violations of the delivery time windows.
2. Penalties for violations of travel time constraints.
3. Cost based on distance travelled.
4. Rewards based on the profit from delivering of the shipment.

The authors put forward several proposals, based on different alternatives for the selection of destinations and the use of the Adaptive Large Neighbourhood Search (ALNS) algorithm as a heuristic:

- Fixed Sample Size (FSS).
- Sample Average Approximation (SAA): $M$ groups are made, each of them having $N$ scenarios. The problem is solved $M$ times. In each iteration, $N$ samples are drawn, independently and identically distributed, from which the difference of that scenario with respect to the optimum is obtained and the average of the differences is calculated. If the deviation and variance are sufficiently small, the solution is accepted.
- Sequential Sample Procedure (SSP): this method is typically used to solve problems though ever-increasing samples. For each iteration, a sample size is selected and observations are chosen to assess the quality of the current solution, based on a statistical estimate of the optimality of the solution. If it is close enough, the method stops; otherwise, the sample size is increased.

### 3.1.2. Studies Based on Data Derived from Real-World Sources

High quality synthetical dataset are common in the literature, but their correspondence with the real world cannot always be assured. The following studies solve this issue by using some real element as a starting point to perform simulations or to obtain synthetically generated data:

- (Janjevic, Winkenbach, & Merchán, 2019) uses the map and shop delivery information of the city of Sao Paulo (Brazil) to carry out their study and alternative proposal.
- (Firdausiyah, Taniguchi, & Qureshi, 2019) use information from the Joint Delivery Systems (JDS) and UCC of the city of Yokohama (Japan) to test the performance of their proposal.

- (Chatterjee, Greulich, & Edelkamp, 2016) use the city of Bremen (Germany) as a model to propose a multimodal shipment distribution model using public transport and electric cargo bikes.
- (Fikar, 2018) uses the locations and means of distribution of a company selling fresh produce (fruit and vegetables) to carry out their simulations and test the advantages of their proposed alternative.
- (Rosano, Demartini, Lamberti, & Perboli, 2018) use information from a distribution company operating in the city of Turin (Italy) to propose an alternative.
- (Chen & Chankov, 2017) use the UPS network of shops in Washington DC (USA) as a basis for testing the performance of their crowdsourced shipping proposal.
- (Muñoz-Villamizar, Montoya-Torres, & Vega-Mejía, 2015) use a group of shops in the city of Bogotá (Colombia) as a model for their proposal.
- (Arslan A., Agatz, Kroon, & Zuidwijk, 2016) uses the UCCs in Rotterdam (Netherlands) to study the impact of the UCCs on the delivery of goods.

## 3.2  Classification of the Bibliography

The state of the art can also be classified according to the technical algorithms used in their solution. Mathematical and heuristic solutions are the predominant choice, or a mixture of both. Their main difference is that, while the mathematical approaches try to find an exact solution by minimising a certain cost function, the heuristic solutions try to reach intermediate, short-term goals. Some of the most notable studies for each case are described sections below

Table 33 describes the basic technical characteristics of the most relevant works found in the state of the art, and the following subsections describe the methods based on their technical nature.

**Table** 3.

| | Static problem | Dynamic problem | Agents' use | Simulation's use | Heuristics' use | Machine learning | Generic algorithms | Mixed-integer linear program. | Non-linear programming |
|---|---|---|---|---|---|---|---|---|---|
| (Janjevic, Winkenbach, & Merchán, 2019) | | | | | ✗ | | | | ✗ |
| (Firdausiyah, Taniguchi, & Qureshi, 2019) | | | ✗ | ✗ | | ✗ | | | ✗ |
| (Rizk & Awad, 2019) | ✗ | | ✗ | | | ✗ | ✗ | ✗ | |
| (Chatterjee, Greulich, & Edelkamp, 2016) | | | ✗ | ✗ | | | | | |
| (Habault, Taniguchi, & Yamanaka, 2018) | | | | | | ✗ | | | |
| (Murray & Chu, 2015) | | | | | ✗ | | | ✗ | |
| (Zhou, Baldacci, Vigo, & Wang, 2018) | | | | | ✗ | ✗ | ✗ | | |
| (Kafle, Zou, & Lin, 2017) | ✗ | | | | ✗ | | | | ✗ |
| (Fikar, 2018) | | ✗ | ✗ | ✗ | ✗ | | | | |
| (Martins-Turner & Nagel, 2019) | | | ✗ | ✗ | | | | | |
| (Rosano, Demartini, Lamberti, & Perboli, 2018) | | ✗ | | ✗ | | | | | |
| (Pimentel & Alvelos, 2018) | | | | | | | | ✗ | |
| (Archetti, Savelsbergh, & Speranza, 2016) | ✗ | | | | ✗ | | | ✗ | |
| (Dayarian & Savelsbergh, 2017) | ✗ | ✗ | | | ✗ | | | | |
| (Gdowska, Viana, & Pedroso, 2018) | | | | | ✗ | | | ✗ | |
| (Li B. , Krushinsky, Woensel, & Reijers, 2016) | | | | | ✗ | ✗ | ✗ | ✗ | |
| (Sampaio-Oliveira, Savelsbergh, Veelenturf, & Woensel, 2018) | | | | | ✗ | | | ✗ | |
| (Chen, Mes, & Schutten, 2018) | ✗ | | | | ✗ | | | ✗ | |
| (Chen & Chankov, 2017) | | | ✗ | ✗ | | | | | |
| (Muñoz-Villamizar, Montoya-Torres, & Vega-Mejía, 2015) | | | | | ✗ | | | ✗ | |
| (Chami, Manier, Manier, & Fitouri, 2017) | | | | | | ✗ | ✗ | ✗ | |
| (Arslan A. , Agatz, Kroon, & Zuidwijk, 2016) | | ✗ | | | | | | ✗ | |

**Technical characteristics of the proposed solutions of the state of the art.**

### 3.2.1   Methods based on mathematical models

Most of the proposals that use mathematical models for the solution of the problem use mixed integer linear programming, although not all of them. It can also be observed that there are two clear tendencies: using the mathematical model to solve the problem and using it only to check the fit of a heuristic solution, although other solutions were also found.

### A)  PROPOSALS THAT USE A PROGRAMMING MODEL FOR THE SOLUTION

(Muñoz-Villamizar, Montoya-Torres, & Vega-Mejía, 2015) use Multi-Depot Capacitated VRP (MDCVRP) models that are solved by mixed integer linear programming, for the comparison between collaborative and non-collaborative scenarios. The computational complexity of the problem is reduced by solving it in two phases. First, the shipping points are associated with the warehouse that serves them and then the vehicle routing problem is solved. In both cases a mixed integer linear programming model is used.

(Pimentel & Alvelos, 2018) include a mixed integer programming model to solve the problem of using urban transport to ship goods.

### B)  PROPOSALS THAT USE A PROGRAMMING MODEL TO CHECK THE QUALITY OF THE PROPOSED SOLUTION

(Murray & Chu, 2015) proposes the Parallel Drone Scheduling TSP (PDSTSP) which uses mixed integer linear programming, although it is used to compare the goodness of heuristics that are proposed to achieve a practical solution.

(Chen, Mes, & Schutten, 2018) also uses mixed integer linear programming for the Multi-Driver, Multi-Parcel Matching Problem (MDMPMP), a model for crowdsourcing dispatch. But given its computational complexity, they develop two heuristics to obtain a practical solution.

(Sampaio-Oliveira, Savelsbergh, Veelenturf, & Woensel, 2018) uses a Multi-Depot Pickup-and-Delivery with Time Window and Transfers (MDPDPTW-T) model formulated as a mixed integer linear program, although its use is discarded and a heuristic is proposed instead.

(Chami, Manier, Manier, & Fitouri, 2017) proposes the Selective Pickup and Delivery Problem with Time Windows and Paired Demands (SPDPTWPD) which is formulated using mixed linear programming, but only its solutions are used to compare with the solution obtained using a hybrid genetic algorithm.

### C)  OTHER

(Taniguchi, Thompson, & Qureshi, 2018) solves the Pickup & Delivery problem through (pure) crowdsourcing in 2 steps. First, a feasibility check of the assignments between requests and deliverers is performed and then the feasible requests are solved by exact methods (TSP-TWPC) under the assumption that there will not be many stops.

### 3.2.2    Methods based on heuristic models

Solutions based on heuristics to obtain a solution to the problem can be classified into two types: those that use an agent-based solution plus the heuristic and those that use only some heuristic solution.

Multi-agent models help to understand the behaviour of different stakeholders in city logistics. They are often used to estimate the social, economic, financial, environmental and energy impact of implementing policy measures in urban areas (Taniguchi, Thompson, & Qureshi, 2018). These stakeholders will be shippers, transport companies, residents, administrations and perhaps other actors such as UCCs.

### A)    HEURISTIC-ONLY SOLUTIONS

(Murray & Chu, 2015) define two heuristics, one for each delivery case defined: one for the case of direct delivery from the warehouse by drone, and another one for the case of using the delivery truck as a mobile base for the drone:

- The first heuristic starts from obtaining a first approximation (heuristic or not) that solves an initial assignment of all deliveries to the truck. Subsequently, the points that the drone can reach are analysed to know if it would be faster to use the drone for the delivery or not (the objective function is to minimize the driver's time).
- The second heuristic assumes that all deliveries reachable by drone will be done by drone and all others will be done by delivery truck and both problems are solved (by heuristic or not). However, as both vehicles must coincide at points along the way, a reallocation of shipments from the drone to the truck is made if this allows minimizing the delivery time.

(Zhou, Baldacci, Vigo, & Wang, 2018) defines a metaheuristic based on a Hybrid Multi-Population Genetic (HMPG) algorithm using a two-level network. A heuristic generates a first set of solutions to the problem, which are classified into feasible and infeasible groups. These groups are then refined by applying genetic operators. It iterates until a solution with a good performance is obtained.

(Kafle, Zou, & Lin, 2017) initially proposes a solution using a mixed integer nonlinear program, for their crowdsourcing delivery problem, which is eventually decomposed into two: a first winner determination problem and a pickup & delivery problem with soft time window. The proposed heuristic is based on Tabu Search, which iteratively solves both problems simultaneously.

The winner determination problem is formulated as a binary integer program and is solved with a branch-&-bound algorithm that does not consider time, and the truck routing problem, which must consider time planning, is solved with an algorithm based on the Simulated Annealing principle to adjust the nodes within the routes. The initial routes are generated using the Nearest Neighbour algorithm.

(Dayarian & Savelsbergh, 2017) considers two solution variants for their crowdsourcing shipping problem. A static one, where it is known in advance when events occur (order arrival, customer arrival, whether the order will arrive, etc.), and a dynamic one, without prior information, assuming stochastic arrivals of both (online) orders and customers. The solution is based on a Tabu Search heuristic to assign one or more shipments to an occasional driver in a way that minimises:

- non-compliance with the delivery time and
- the transport cost.

The static solution, which knows all the information and can use it to adjust the allocation of consignments to an occasional driver or to its own delivery staff.

The dynamic solution, which does not know the future, is approached in two different ways: without any information about the future and therefore solving the problem known at the moment and extrapolating the current arrival rate (of orders and customers) into the future.

(Li B., Krushinsky, Woensel, & Reijers, 2016) proposes a problem which consists of maximizing the transport of passengers and parcels by cab, for which a two-stage stochastic programming model is proposed. To solve it, an Adaptive Large Neighbourhood Search (ALNS) heuristic is proposed together with 3 strategies for scenario generation. The new solution is accepted or not according to a Simulated Annealing process.

(Sampaio-Oliveira, Savelsbergh, Veelenturf, & Woensel, 2018) proposes a solution to the multi-warehouse problem with pickup-&-delivery with time windows and transfers using the ALNS heuristic which allows to improve the performance of the exact solution. To control the acceptance of a new solution, they rely on the Simulated Annealing process, although sometimes worse solutions are accepted, in order to increase the diversity of options.

(Chen, Mes, & Schutten, 2018) proposes to use the daily commuting of the population (e.g., to work) to reduce the cost of delivery through a multi-parcel multi-driver assignment, using transfers. It is first formulated using integer programming, but also using two heuristics in view of the complexity of the solution:

- The first heuristic, the time-compatibility heuristic (TC-heuristic) assigns each dispatch to the calculated shortest path (an A* algorithm is used to obtain the shortest path and the alternatives that involve less than a certain detour) that allows reaching the destination on time. This algorithm can be too costly.
- The second heuristic, called time-expanded graph based (TEG-heuristic) aims to be able to solve problems of a more realistic size. The driver availability information is used to construct the TEG-graph and then the drivers are assigned to the dispatches.

A key point in both cases is the compatibility that can exist between the solutions calculated by the A* algorithm and the constraint of always keeping the dispatch in possession of one of the drivers, which implies taking into account waits in the transfers of the dispatches.

(Chami, Manier, Manier, & Fitouri, 2017) proposes a solution to a variant of the Selective Pickup and Delivery Problem (PDP) with Time Windows and Paired Demands (SPDPTWPD) problem which is tackled by means of a hybrid genetic heuristic. In this case, the problem also includes constraints such as maximum capacity of delivery vehicles due to access zone restrictions. The proposed solution uses a hybrid genetic algorithm with a local search that helps to minimise distance and maximise economic performance. The mutation operator is based on a local search.

### B)  AGENT MODEL PLUS HEURISTIC SOLUTIONS

(Chatterjee, Greulich, & Edelkamp, 2016) proposes a solution based on the use of existing public transportation networks to achieve benefits such as minimizing congestion or reducing greenhouse gases. Cost, time and energy are optimized using a multi-agent simulation model. Dijkstra's minimum path algorithm and a Nested Monte Carlo Search are used to construct the cost and time matrix required to generate the intermodal route.

(Martins-Turner & Nagel, 2019) approach is based on finding out the effect of being able to make more than one trip with the same vehicle, within its schedule. The proposed solution is applied to two different scenarios, with different characteristics: food delivery and parcel delivery. In their Vehicle Routing Problem (VRP) two different problems are identified: (1) assigning requests to routes and (2) determining the optimal sequence in which these requests are served within the route (routing).

MATSim (Multi-Agent Transport Simulation) is used to solve the problem, which iteratively performs three steps: traffic simulation, score calculation and replanning. Agents performing activities are simulated. Once the simulation is done, the agents try to improve their plans. The simulation process stops once a limit of iterations is reached. Internally, a VRP is solved to perform routing.

(Rizk & Awad, 2019) proposes to solve their Pickup and Delivery Problem with Coalition Formation (PDP-CF). They formulate the mixed integer programming model, although its solution is problematic because the search space grows exponentially with the number of vehicles. To solve it, they propose a quantum genetic algorithm whose tuning function is based on the objective function of the linear program.

(Fikar, 2018) proposes a solution for the delivery of food (fruit and vegetables) by means of a professional fleet of delivery drivers includes dynamic variations in the route. If new orders appear and can be included in the already calculated route. The objective is to minimise the route travelled and maximise the shelf life of perishable goods. The problem consists of the following decisions:

- choose the point that serves the order.
- choose the products.
- choose the vehicle.

- deciding the route.
- the time to reach the destination.

The solution is based on applying geographic information with simulation and optimisation procedures. Demand is generated randomly, and decisions are made dynamically using a heuristic optimisation procedure. As requests arrive, the heuristic tries to insert them into an already running route, if possible. The process is iterative, adding new requests and replanning routes until there is no improvement.

### 3.2.3    Methods based on a mixture of mathematical and heuristics models

Several authors have proposed to take advantages of both of the previous approaches, by developing solutions based on a mixture of mathematical and heuristic methods:

(Firdausiyah, Taniguchi, & Qureshi, 2019) proposes an agent- and simulation-based model to adjust their pricing model for shipping services. However, to evaluate the cost functions, they use the VRPSTW (Vehicle Routing Problem with Soft-Time Window) model proposed in (Ali Gul Qureshi, 2010). The authors propose to use an agent-based model and simulations to model the interaction between the UCC and various transportation companies in order to adjust service prices. To calculate costs, the VRPSTW model is used, while the simulation process allows the application of learning techniques (Reinforcement Learning) to make the price adjustment and allow adaptation to changing conditions.

(Archetti, Savelsbergh, & Speranza, 2016) replaces the initially presented integer program by a multi-start heuristic combining neighbourhood and tabu searches and several simpler integer programs. The heuristic part of the solution uses several (5) initial solutions, where the first (a) is the solution that does not use occasional drivers, which adds drivers in increasing order of distance to the store. Then integer programs are used to (b) determine which customers are served by occasional drivers, (c) maximize the demand served by occasional drivers, (d) determine the maximum number of occasional drivers, and (e) to minimize the cost. Tabu search selects dispatches to reroute or to assign to occasional drivers to try to improve the solution. For the occasional drivers that have no assignment, another integer program is used that creates a new assignment, which allows to modify the solution reached.

(Gdowska, Viana, & Pedroso, 2018) first applies a deterministic solution by solving the problem with the VRPOD (Vehicle Routing Problem with Occasional Drivers) model using integer programming to have a reference and then applies a second model to minimize the cost by assigning shipments to occasional drivers, the cost function optimizes taking into account the randomness in the response of occasional drivers.

(Muñoz-Villamizar, Montoya-Torres, & Vega-Mejía, 2015) performs a comparison between the collaboration and non-collaboration scenarios (between several retail delivery companies). Moreover, the use of integer programming is proposed, specifically, a Capacitated Vehicle Routing Problem (CVRP) for the non-collaboration scenario and a Multi-Depot Capacitated Vehicle Routing Problem (MDCVRP) for the collaboration scenario. To reduce the complexity of the MDCVRP problem, the use of two heuristic functions in hierarchy is proposed. First, the dispatch points to be served from each warehouse are assigned and in a second step, the routing problem is solved.

(Janjevic, Winkenbach, & Merchán, 2019) The solution to delivery using a multi-level network and CDPs, based on nonlinear models, which are responsible for determining:

- whether satellites are active.
- whether a CDP is included in a route.
- whether a route is served from a satellite or not, becomes very expensive when there are many dispatches.

To avoid this, a heuristic solution is proposed based on:

- first setting the configuration for CDPs.
- second setting the configuration for the satellites
- then solving the routing problem for the shipments assuming known satellite configurations, and pickup and delivery points, which allows to reduce its nonlinear programming model to a linear programming one.

The final algorithm acts through two searches, one breath-first and one deep-first, to fix the satellite and CDP configurations and then applies the linear program to evaluate the cost, keeping the configurations fixed.

### 3.2.4    Static vs Dynamic Solutions

Most of the solutions proposed in the references work with the static nature of the problem: the shipments, the delivery drivers and the intervals for deliveries (if any) are fully determined when the workload allocation is made and the route is calculated.

On the other hand, some authors prefer to focus on the dynamic nature of the problem: orders or delivery drivers (or both) arrive randomly, and it is not known what or when will arrive in advance. As noted in (Gdowska, Viana, & Pedroso, 2018), crowdsourcing-based solution proposals are likely easier to handle by dynamic methods given their nature.

(Fikar, 2018) describes a fresh produce (fruits and vegetables) shipping problem that includes a quality loss model based on time and temperature of the product, so it has tighter temporal constraints than other problems. In this case, a decision support system plans the routes of professional delivery drivers - several delivery drivers and several shops as a point of origin - to, among other things, limit waste and maintain the quality of the delivered product. The dynamics are incorporated in two ways:

1. Making schedules with the information available up to a certain point in time, without knowing the future.
2. Trying to include new requests in routes that are already formed: an already calculated vehicle route can be modified to accommodate a new order (going to a shop to pick it up and deliver it to the destination) if appropriate.

(Rosano, Demartini, Lamberti, & Perboli, 2018) describes a crowdsourcing problem. Their proposal is based on a platform and a mobile application (UFreight) targeting B2B (Business-to-Business) and B2C (Business-to-Client) segments. New requests are notified to the most compatible shipping companies (by geographic location and capacity) that bid for the shipment. The customer making the request can choose the offer from among those presented. They use a simulation-optimisation environment, described in (Perboli, Rosano, Saint-Guillain, & Rizzo, 2018), to solve the Dynamic and Stochastic Vehicle Routing Problem with Time Windows (DS-VRPTW). The DS-VRPTW model is described in (Saint-Guillain, Deville, & Solnon, 2015) and allows modelling the existence of random requests that can be included in routes that have already been initiated. At each simulation-optimisation step, the action to be taken is chosen based on known requests and a known probability of occurrence of requests.

In (Dayarian & Savelsbergh, 2017) a solution for crowdsourcing delivery based on occasional drivers (customers of a shop) is proposed. In this case, both the requests received online and the customers of the shop will have random characteristics: it is not known how many will appear (neither requests nor customers) and when, in principle, although some assumptions can be made: one can assume that the current ratios are maintained in the future, for example. The authors compare the result of assigning existing requests to customers as they arrive (with no more information than is available at the time), with the possibility of incorporating data on the distribution of request and customer arrivals, and conclude that the latter improves the performance of the system. In any case, in their model, the authors include a backup system for the delivery system: the shop has dedicated deliverymen who take care of the requests that cannot be assigned to anyone else.

(Arslan A., Agatz, Kroon, & Zuidwijk, 2016) also includes a crowdsourcing solution to the Pickup & Delivery problem, also with drivers on their daily commute. A platform automatically creates the allocation between deliveries and drivers. In this case, it is considered that there may be several stops on the route, for pickups or deliveries, so a routing solution is needed in addition to the driver assignment. The backup solution proposed to be used is an external delivery company, so there is no fleet of its own. Both requests and drivers will arrive randomly, so a solution that allows for this dynamic in the process will be necessary. To this end, the authors propose to use the arrival times of requests and customers to launch a process that determines the feasibility of a shipment with the customer routes (and perhaps the routing, if there are feasible shipments) of the orders pending shipment. The authors study two different cases, depending on the limits of the delivery interval. The first case is to make the assignment as soon as possible, while the second case is to wait until the upper limit of the interval is approached, in the hope that more than one shipment can be assigned to the same driver, allowing better assignments to be made.

## 4. Overall Categorization

In (Braekers, Ramaekers, & Nieuwenhuyse, 2016) a review of the literature for the Vehicle Routing Problem (VRP) up to the date of publication is carried out, and a taxonomy is created for them, based on five basic characteristics: study type, scenario characteristics, physical problem characteristics, information characteristics and data characteristics.

Taking inspiration from their taxonomy, a classification is proposed. Table 4. Proposed classification of the literature.4 summarizes all the literature solutions and categorizes them, grouping them according to their approach method and their specific solution.

| Approach | Solution | Bibliographic works |
|---|---|---|
| Deliverymen | Hired deliverymen | (Janjevic, Winkenbach, & Merchán, 2019), (Firdausiyah, Taniguchi, & Qureshi, 2019), (Chatterjee, Greulich, & Edelkamp, 2016), (Murray & Chu, 2015), (Fikar, 2018), (Martins-Turner & Nagel, 2019), (Rosano, Demartini, Lamberti, & Perboli, 2018), (Muñoz-Villamizar, Montoya-Torres, & Vega-Mejía, 2015) |
| | Crowdsourcing / Occasional drivers | (Habault, Taniguchi, & Yamanaka, 2018), (Kafle, Zou, & Lin, 2017), (Sampaio-Oliveira, Savelsbergh, Veelenturf, & Woensel, 2018), (Chen & Chankov, 2017), (Arslan A. , Agatz, Kroon, & Zuidwijk, 2016) |
| | Mixture of hired deliverymen and crowdsourcing | (Archetti, Savelsbergh, & Speranza, 2016), (Dayarian & Savelsbergh, 2017), (Gdowska, Viana, & Pedroso, 2018), (Chen, Mes, & Schutten, 2018) |
| | Other | (Pimentel & Alvelos, 2018), (Li B. , Krushinsky, Woensel, & Reijers, 2016) |
| Temporal restrictions | Few hours (1-2 h) | (Habault, Taniguchi, & Yamanaka, 2018), (Fikar, 2018), (Archetti, Savelsbergh, & Speranza, 2016), (Li B. , Krushinsky, Woensel, & Reijers, 2016) |
| | All day long | (Janjevic, Winkenbach, & Merchán, 2019), (Firdausiyah, Taniguchi, & Qureshi, 2019), (Chatterjee, Greulich, & Edelkamp, 2016), (Murray & Chu, 2015), (Zhou, Baldacci, Vigo, & Wang, 2018), (Kafle, Zou, & Lin, 2017), (Martins-Turner & Nagel, 2019), (Rosano, Demartini, Lamberti, & Perboli, 2018), (Pimentel & Alvelos, 2018), (Archetti, Savelsbergh, & Speranza, 2016), (Dayarian & Savelsbergh, 2017), (Sampaio-Oliveira, Savelsbergh, Veelenturf, & Woensel, 2018), (Chen, Mes, & Schutten, 2018), (Chen & Chankov, 2017), (Muñoz-Villamizar, Montoya-Torres, & Vega-Mejía, 2015), (Chami, Manier, Manier, & Fitouri, 2017), (Arslan A. , Agatz, Kroon, & Zuidwijk, 2016) |
| Additional cost-saving mechanisms | No mechanism | (Rizk & Awad, 2019),  (Murray & Chu, 2015), (Fikar, 2018), (Rosano, Demartini, Lamberti, & Perboli, 2018), (Archetti, Savelsbergh, & Speranza, 2016), (Dayarian & Savelsbergh, 2017), (Gdowska, Viana, & Pedroso, 2018), (Chen & Chankov, 2017), (Muñoz-Villamizar, Montoya-Torres, & Vega-Mejía, 2015), (Chami, Manier, Manier, & Fitouri, 2017), (Arslan A. , Agatz, Kroon, & Zuidwijk, 2016) |
| | UCCs | (Firdausiyah, Taniguchi, & Qureshi, 2019), (Martins-Turner & Nagel, 2019) |
| | Echelon network, CDPs | (Janjevic, Winkenbach, & Merchán, 2019), (Zhou, Baldacci, Vigo, & Wang, 2018), (Kafle, Zou, & Lin, 2017), (Martins-Turner & Nagel, 2019), (Pimentel & Alvelos, 2018) |
| | Alternative options (multimodal transport, colaboration, transfers, etc) | (Chatterjee, Greulich, & Edelkamp, 2016), (Habault, Taniguchi, & Yamanaka, 2018), (Muñoz-Villamizar, Montoya-Torres, & Vega-Mejía, 2015), (Sampaio-Oliveira, Savelsbergh, Veelenturf, & Woensel, 2018), (Chen, Mes, & Schutten, 2018) |
| Destinations | One or two stops | (Rizk & Awad, 2019), (Habault, Taniguchi, & Yamanaka, 2018), (Dayarian & Savelsbergh, 2017), (Gdowska, Viana, & Pedroso, 2018), (Sampaio-Oliveira, Savelsbergh, Veelenturf, & Woensel, 2018), (Chen & Chankov, 2017) |
| | Few stops | (Firdausiyah, Taniguchi, & Qureshi, 2019), (Kafle, Zou, & Lin, 2017), (Fikar, 2018), (Pimentel & Alvelos, 2018), (Archetti, Savelsbergh, & Speranza, 2016), (Li B. , Krushinsky, Woensel, & Reijers, 2016), (Chen, Mes, & Schutten, 2018), (Arslan A. , Agatz, Kroon, & Zuidwijk, 2016) |
| | Route with many stops | (Janjevic, Winkenbach, & Merchán, 2019), (Chatterjee, Greulich, & Edelkamp, 2016), (Murray & Chu, 2015), (Zhou, Baldacci, Vigo, & Wang, 2018), (Martins-Turner & Nagel, 2019), (Rosano, Demartini, Lamberti, & Perboli, 2018), (Muñoz-Villamizar, Montoya-Torres, & Vega-Mejía, 2015), (Chami, Manier, Manier, & Fitouri, 2017) |
| **Vehicles** | Vans, lorries, … | (Janjevic, Winkenbach, & Merchán, 2019), (Firdausiyah, Taniguchi, & Qureshi, 2019), (Zhou, Baldacci, Vigo, & Wang, 2018), (Fikar, 2018), (Martins-Turner & Nagel, 2019), (Rosano, Demartini, Lamberti, & Perboli, 2018), (Muñoz-Villamizar, Montoya-Torres, & Vega-Mejía, 2015), (Chami, Manier, Manier, & Fitouri, 2017) |
| | Multi-modal (public transport) | (Chatterjee, Greulich, & Edelkamp, 2016), (Pimentel & Alvelos, 2018), (Li B. , Krushinsky, Woensel, & Reijers, 2016) |
| | Personal vehicles, motocycle, bicycles, walking | (Habault, Taniguchi, & Yamanaka, 2018), (Kafle, Zou, & Lin, 2017), (Martins-Turner & Nagel, 2019), (Archetti, Savelsbergh, & Speranza, 2016), (Dayarian & Savelsbergh, 2017), (Gdowska, Viana, & Pedroso, 2018), (Sampaio-Oliveira, Savelsbergh, Veelenturf, & Woensel, 2018), (Chen, Mes, & Schutten, 2018), (Chen & Chankov, 2017), (Arslan A. , Agatz, Kroon, & Zuidwijk, 2016) |
| | UAV (Drones) and similar | (Murray & Chu, 2015), (Rizk & Awad, 2019) |
| **Solution type** | Mathematical | (Pimentel & Alvelos, 2018), (Muñoz-Villamizar, Montoya-Torres, & Vega-Mejía, 2015) |
| | Mixture of mathematics and heuristical | (Janjevic, Winkenbach, & Merchán, 2019), (Firdausiyah, Taniguchi, & Qureshi, 2019), (Archetti, Savelsbergh, & Speranza, 2016), (Gdowska, Viana, & Pedroso, 2018), (Muñoz-Villamizar, Montoya-Torres, & Vega-Mejía, 2015) |
| | Heuristic | (Rizk & Awad, 2019), (Chatterjee, Greulich, & Edelkamp, 2016), (Habault, Taniguchi, & Yamanaka, 2018), (Murray & Chu, 2015), (Chen, Mes, & Schutten, 2018), (Archetti, Savelsbergh, & Speranza, 2016), (Sampaio-Oliveira, Savelsbergh, Veelenturf, & Woensel, 2018), (Chami, Manier, Manier, & Fitouri, 2017), (Zhou, Baldacci, Vigo, & Wang, 2018), (Kafle, Zou, & Lin, 2017), (Dayarian & Savelsbergh, 2017), (Li B. , Krushinsky, Woensel, & Reijers, 2016), (Taniguchi, Thompson, & Qureshi, 2018), (Chatterjee, Greulich, & Edelkamp, 2016), (Martins-Turner & Nagel, 2019), (Fikar, 2018), (Chen & Chankov, 2017) |
| **Problem type** | Delivery | (Janjevic, Winkenbach, & Merchán, 2019), (Firdausiyah, Taniguchi, & Qureshi, 2019), (Chatterjee, Greulich, & Edelkamp, 2016), (Habault, Taniguchi, & Yamanaka, 2018), (Zhou, Baldacci, Vigo, & Wang, 2018), (Kafle, Zou, & Lin, 2017), (Fikar, 2018), (Martins-Turner & Nagel, 2019), (Pimentel & Alvelos, 2018), (Archetti, Savelsbergh, & Speranza, 2016), (Dayarian & Savelsbergh, 2017), (Gdowska, Viana, & Pedroso, 2018), (Li B. , Krushinsky, Woensel, & Reijers, 2016), (Chen, Mes, & Schutten, 2018), (Chen & Chankov, 2017), (Muñoz-Villamizar, Montoya-Torres, & Vega-Mejía, 2015), (Arslan A. , Agatz, Kroon, & Zuidwijk, 2016) |
| | Delivery and collection | (Rizk & Awad, 2019), (Rizk & Awad, 2019), (Rosano, Demartini, Lamberti, & Perboli, 2018), (Sampaio-Oliveira, Savelsbergh, Veelenturf, & Woensel, 2018), (Chami, Manier, Manier, & Fitouri, 2017) |
| **Sent product type** | Perishable (Food) | (Habault, Taniguchi, & Yamanaka, 2018), (Fikar, 2018), (Martins-Turner & Nagel, 2019) |
| | Other (packages, letter, …) | (Janjevic, Winkenbach, & Merchán, 2019), (Firdausiyah, Taniguchi, & Qureshi, 2019), (Chatterjee, Greulich, & Edelkamp, 2016), (Murray & Chu, 2015), (Zhou, Baldacci, Vigo, & Wang, 2018), (Kafle, Zou, & Lin, 2017), (Martins-Turner & Nagel, 2019), (Rosano, Demartini, Lamberti, & Perboli, 2018), (Pimentel & Alvelos, |

| | | 2018), (Archetti, Savelsbergh, & Speranza, 2016), (Dayarian & Savelsbergh, 2017), (Gdowska, Viana, & Pedroso, 2018), (Li B. , Krushinsky, Woensel, & Reijers, 2016), (Sampaio-Oliveira, Savelsbergh, Veelenturf, & Woensel, 2018), (Chen, Mes, & Schutten, 2018), (Chen & Chankov, 2017), (Muñoz-Villamizar, Montoya-Torres, & Vega-Mejía, 2015), (Chami, Manier, Manier, & Fitouri, 2017), (Arslan A. , Agatz, Kroon, & Zuidwijk, 2016) |
|---|---|---|

**Table 4. Proposed classification of the literature.**

## 5. Urban Road Freight Impacts

Another important aspect to consider when analysing the Last Mile Delivery problem is the impact of urban freight and the existing literature which implements solutions for reducing its impact.

(Coulombel N. *et al.*, 2018) proposes an environmental impact study of urban road freight focused on the pollutant emissions. It takes place in the Paris region. They develop a freight demand model, a multiclass traffic assignment model and a road emission model. They value the environmental impact of urban road freight in approximately €2.1 billion, per year.

One possible solution is electric mobility. (Iwan, S., *et al.*, 2019) analyze the situation and projects in the European market. Recent technological developments have made these vehicles suitable for city logistics, with the advantage of reduced exhaust and noise emissions. The main drawback at the moment was that logistic players are financially reluctant to make the change due to the uncertain depreciation rate of electric freight vehicles. However, the situation showed some improvement over the course of the own project.

(İmre, Ş., *et al.*, 2021) focuses on addressing the stakeholder needs regarding Electric Vehicles (EVs) in urban freight transport. The study is focused on Türkiye. They recognise a lack of knowledge on specific technical and organizational issues, which causes the commercial fleets to be far below their potential. They conclude that EVs adoption is limited by a lack of easily available and reliable information, generating uncertainties and perceived risks. Other factors are also identified.

Moreover, the integration of new technologies into city logistics is also considered in the literature. (Comi, A. & Russo, F., 2022) study the current iteration of cities with the Internet of Things (IoT), block chain, big data, and artificial intelligence, which in turn generate integrated and dynamic city logistics solutions. They assert that routing and scheduling can benefit from IoT for daily decision making, and from big data for weekly decision making. Moreover, blockchain allows to monitor the provided service with high trust, and to monitor preservation conditions of perishable goods. A closer collaboration of the public administrations and freight operators is proposed as it would benefit all stakeholders.

## 6. LMD Proposal

Based on the analysed works of the state of the art, we propose a Multi-Agent System (MAS) with a mixed system of heuristics and mathematical methods for solving the vehicle routing problem. Its main goal is to optimise the delivery cost of non- perishable products, and to have a good scalability. This configuration is chosen as the considered problem is unstable, with products arriving to warehouses in an unpredictable order, and as a good scalability is desired to expand the delivery business if necessary. The preferred delivery vehicle is a van, and trucks can optionally be used too.

Packages will be sent from the distribution centre to smaller storage centres (local businesses) using trucks or vans, forming clusters which centroid is the storage/distribution centre. Let $a_i$ be the distribution point, where $a_0$ is the centroid and $N$ is the number of distribution points. Then the goal is to find:

$$\min\left(\sum_{i=0}^{N} d(a_i, a_{i+1})\right) \tag{1}$$

Let $c_j$ be the j-th centroid, with $j \in \mathbb{R}^M$ and $M$ being the total number of centroids. Then, package $p_k$ with destination $a_{dest}$ will be stored in the centroid j-th such as:

$$\min_j\left(d(c_j, a_{dest})\right) \tag{2}$$

Packages will be stored in such centroid for up to 2 days. The recipient of the package then has the possibility of collecting the parcel within these two days, or to wait until the third day, when it will be delivered to his chosen address.

Every day, the distribution route is calculated before the start of the working day. In this route, vans collect all the third-day packages in the centroid, distributes them to their destinations, then move to the next centroid and repeat the process. Vans' first centroid of the day is $o_x$ such as:

$$\max\left(\sum_i d(o_i, o_{i+1})\right) \tag{3}$$

for $x \in \{1, \dots, P\}$ with $P$ being the number of available vans. The path from the vans' depot to their respective initial centroid is calculated using an adapted version of the Dijkstra's algorithm for finding the shortest paths between nodes in a graph:

1. Let $v_i$ be the vertices/nodes and let $w(v_i, v_j)$ be the estimated travel time between adjacent nodes $v_i$ and $v_j$.
2. Let's initialize $d(v_0, v_i) = 0$, where $d(v_0, \cdot)$ is an array containing the estimated travel times between the origin and any node. Let's also initialize $d(v_0, v_j) = \infty$, for all $i, j \neq 0$.
3. Let $Q$ be a set containing all the nodes in the graph, and let $S$ be an empty set.
4. While $Q$ is not an empty set, follow the process below:

a. Select the node $v \in Q$ verifying that $\min_i(d(v_0, \ v))$ and that $v \notin S$.

b. Add the node $v$ to $S$.

c. Update the value of $d(v_0, \ v_i)$, for all $i$, as $d(v_0, \ v) + weight(v, \ v_i)$ if $d(v_0, \ v) + weight(v, \ v_i) < d(v_0, \ v_i)$.

When the algorithm has visited all the nodes, the smallest distance to the destination is found. Vans will follow this route to the origin and stop in the centroids (nodes) of their path, making their corresponding deliveries. Incidentally, if a centroid has been previously visited during the day by another van, it will be ignored as no deliveries as left.

The order of the centroids is calculated based on the nearest neighbour (NN) algorithm: all the centroids with third-day packages must be visited by the delivery vehicle by the end of the day, so the initial path contains these centroids as nodes. If the estimated overall time is calculated to be lower than the worker's shift, then the closest nodes are added to the graph and the solution re-calculated; until no additional nodes can be added.

The "neighbourhood route" (route from one centroid its destinations) is planned using the meta-heuristics proposed in (Hasan, M., & Niyogi, R.; 2020), adding an extra objective to its cost function: it must take into account that the end of the route should be as close as possible to the following centroid. Therefore, the objective function is:

$$\min \left( \sum_{k \in V_m} \sum_{j \in D \cup C} \sum_{i \in D \cup C} FC_k * x_{ij}^{mk} + \sum_{k \in V_m} \sum_{j \in D \cup C} \sum_{i \in D \cup C} d_{ij} * x_{ij}^{mk} * VarC_k \right) \tag{4}$$

Where $FC_k$ is the fixed cost of vehicle $k$, $x_{ij}^{mk}$ is a decision variable, $VarC_k$ is the variable cost of vehicle $k$, and $d_{ij}$ is the distance between customer $i$ and $j$. If $i, j \in \{1, \dots, N\}$; then the starting and end points will be adding to the set, resulting in $i, j \in \{0, \dots, N + 1\}$. All the constrains and penalties' details can be found in (Hasan, M., & Niyogi, R.; 2020).

The failed deliveries will be unloaded in the following centroid, and a new delivery attempt will be performed on the next day. The order of the centroids is reversed every day in order to allow this. Regarding the neighbourhood route planning, the mis-placed packages will be ignored when calculating the route of their new location, and they will be added to the calculation of the route of their original location.

It is also important to plan for anomalies in the number of packages as they are common in real-world applications. In case the number of packages is higher than usual on a particular day, some of the received packages should be delivered on their second and first day. If the average daily number of packages increases, more vehicles should be used.

This former scenario is where the MAS is applied: each new vehicle will be a new agent which will communicate with other agents to avoid duplicities. At the start of the day, the centroids which maximise the initial distance between agents are identified, and they will start their shift by then (avoiding overlapping issues). Figure 1 shows a high-level representation of the system. The main benefits of using a MAS system are: good scalability, high reliability, high flexibility, and the high degree of autonomy of each vehicle.

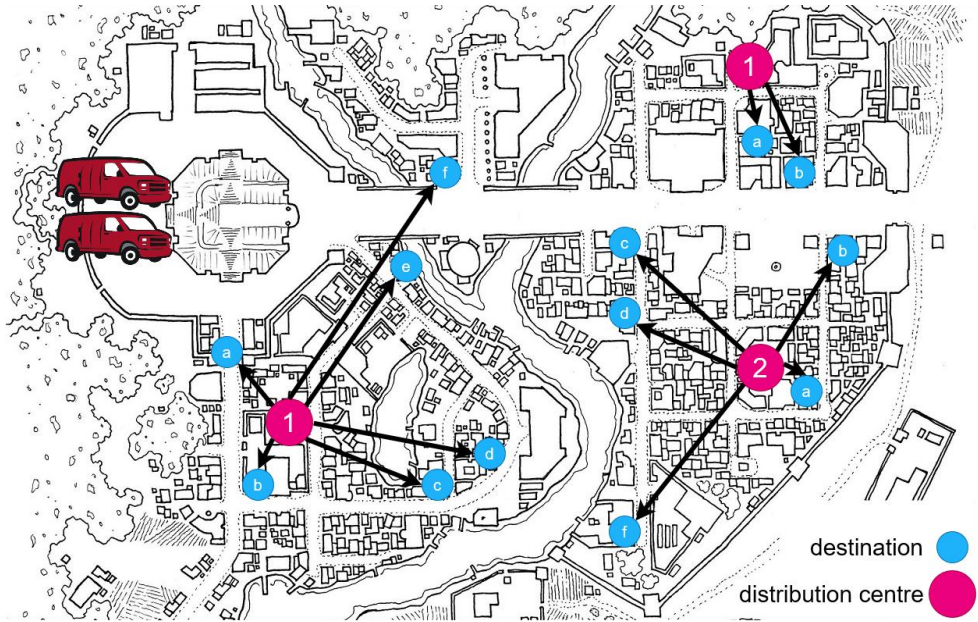

**Figure 1. Representation of the proposed system. Distribution centres and destinations are numbered by agent's visit order, and black lines represent the destinations of the packages they stored.**

The proposed system encourages users to collect the packages themselves, reducing the need of individual deliveries to the distributor. Moreover, the number of journeys to the central storage centre is reduced, which minimizes the petrol consumption and maximizes the expected per-hour-productivity of the workers.

The main challenges of its implementations are finding local shops which accept to be used as storage centres for an economical compensation; and maintaining good customer satisfaction levels as the drawback of cost optimization is a longer delivery time in some cases.

## 7. Conclusions

Last mile solutions make it possible for logistic actors to make decisions regarding parcel deliveries aiming to optimise the financial and environmental costs. As the studied state-of-the-art has shown, participants using concepts of last mile delivery achieved considerable economic savings. Urban road traffic is considered as a particular case, where EVs have a great potential yet to be exploited.

The article is addressed to the service providers, manufacturing actors and sellers with online presence, proposing a novel classification of the bibliography (Section 4) and a description of the most prominent works according to certain technical aspects. Different algorithmical approaches are needed for different scenarios, with no generalist model which adapts well to all situations. More attention to this research question is needed due to its high impact on the modern world, as it is under-developed by comparison to other major parts of the current supply chain.

This research work aims to ease the development of studies of new LMD works. The literature is classified according to some of its most fundamental characteristics, challenges of the field are identified and a theoretical LMD proposal is put forward. All in all, we provide a comprehensive summary of the last mile delivery field which contributes to a more coherent understanding of the topic.

Data sources (real vs synthetical), heuristic/mathematical algorithms and static/dynamic solutions have been identified as the most defining characteristics of the works in the literature, with each of the choices being suitable for specific situations.

A future research line would consist of expanding this research classification, creating a subdivision of each category and increasing the number of categories. Moreover, we would like to study the performance of the LMD proposal in the real world and analyse which modifications would be needed to adapt it to new situations. Overall, the long term goal is to understand to key characteristics of how models adapt to new situations, and put forward a more generalist system, which could be widely used among the industry with little modifications.

## 8. ACKNOWLEDGMENTS

This research has been supported by the project "Intelligent and sustainable mobility supported by multi-agent systems and edge computing (InEDGE-Mobility): Towards Sustainable Intelligent Mobility: Blockchain-based framework for IoT Security", Reference: RTI2018–095390-B-C32, financed by the Spanish Ministry of Science, Innovation and Universities (MCIU), the State Research Agency (AEI) and the European Regional Development Fund (FEDER).